\documentclass[a4paper]{article}
\usepackage[dvips]{color,graphicx}
\usepackage{psfrag}
\usepackage{amsmath,amssymb,amsthm,latexsym}

\newcounter{theocounter}
\newtheorem{Thm}[theocounter]{Theorem}
\newtheorem{Prop}{Proposition}[section]
\newtheorem{Cor}[Prop]{Corollary}
\newtheorem{Lemma}[Prop]{Lemma}
\newtheorem{Definition}[Prop]{Definition}
\newtheorem*{Thm*}{Theorem}
\newtheorem*{Remark}{Remark}

 \newlength\headseptemp

\newcommand{\Lp}{\widetilde{\Delta}}
\newcommand{\alg}{\widetilde{\al}}

\newcommand{\dd}{{\partial}}

\newcommand{\al}{{\alpha}}

\newcommand{\ka}{{\kappa}}
\newcommand{\curv}{{\kappa}}
\newcommand{\ph}{{\varphi}}

\newcommand{\lm}{{\lambda}}

\newcommand{\ess}{{\mathrm {ess}}}
\newcommand{\lra}{{\longrightarrow}}

\newcommand{\R}{{\mathbb R}}

\newcommand{\Z}{{\mathbb Z}}

\newcommand{\V}{{ V}}
\newcommand{\E}{{ E}}
\newcommand{\F}{{ F}}
\newcommand{\G}{{ G}}
\newcommand{\T}{{ T}}
\newcommand{\C}{{ C}}
\newcommand{\Cut}{{\rm Cut}}
\newcommand{\vol}{{\rm vol}}

\newcommand{\Hm}[1]{\leavevmode{\marginpar{\tiny%
$\hbox to 0mm{\hspace*{-0.5mm}$\leftarrow$\hss}%
\vcenter{\vrule depth 0.1mm height 0.1mm width \the\marginparwidth}%
\hbox to 0mm{\hss$\rightarrow$\hspace*{-0.5mm}}$\\\relax\raggedright
#1}}}

\date{\today}

\title{Cheeger constants, growth and spectrum of locally tessellating
  planar graphs}
\author{Matthias Keller \thanks{e-mail:
    m.keller@uni-jena.de} \\
  Mathematical Institute\\
  FSU Jena\\
  D-07743 Jena, Germany
  \and
  Norbert Peyerimhoff \thanks{ e-mail: norbert.peyerimhoff@durham.ac.uk} \\
  Department of Math. Sciences\\
  University of Durham\\
  Durham DH1 3LE, UK}
\begin{document}

\maketitle


\begin{abstract} \noindent In this article, we study relations between
  the {\em local geometry} of planar graphs (combinatorial curvature)
  and {\em global geometric invariants}, namely the Cheeger constants
  and the exponential growth. We also discuss spectral applications.
\end{abstract}


\section{Introduction}

A {\em locally tessellating planar graph} $\G$ is a tiling of the
plane with all faces to be polygons with finitely or infinitely many
boundary edges. The edges of $\G$ are continuous rectifiable curves
without self-intersections. Faces with infinitely many boundary edges
are called infinigons and occur, e.g., in the case of planar
trees. The sets of vertices, edges and faces of $\G$ are denoted by
$\V, \E$ and $\F$ (see the beginning of Section \ref{sec:proofthm1} for
precise definitions). The function $d(v,w)$ denotes the combinatorial distance
between two vertices $v,w \in \V$, where each edge is assumed to have
combinatorial length one. For any pair $v,w$ of adjacent vertices we
write $v \sim w$.

Useful {\em local concepts} of a planar graph $\G$ are combinatorial
curvature notions. The finest curvature is defined on the corners of
$\G$. A corner is a pair $(v,f) \in \V \times \F$, where $v$ is a
vertex of the face $f$. The {\em corner curvature} $\ka_C$ is defined
as
$$ \ka_C(v,f) = \frac{1}{|v|} + \frac{1}{|f|} - \frac{1}{2}, $$
where $|v|$ and $|f|$ denote the degree of the vertex $v$ and the face
$f$. If $f$ is an infinigon, we set $|f| = \infty$ and $1/|f| =
0$. The {\em curvature at a vertex} $v \in \V$ is given by
$$ \ka(v) = \sum_{f: v \in f}\ka_C(v,f) = 1 - \frac{|v|}{2} +
\sum_{f: v \in f} \frac{1}{|f|}. $$

For a finite set $W\subset \V$ we define $\ka(W) = \sum_{v \in W}
\ka(v)$ and the {\em average vertex curvature} by
$$
\overline{\ka}(W) = \frac{1}{|W|} \ka(W).
$$
These combinatorial curvature definitions arise naturally from
considerations of the Euler characteristic and tessellations of closed
surfaces, and they allow to prove a combinatorial Gau{\ss}-Bonnet
formula (see \cite[Theorem 1.4]{BP1}). Similar combinatorial curvature
notions have been introduced by many other authors, e.g.,
\cite{St,Gro,Woe,Hi}. Let us already mention two global geometric
consequences of the curvature sign:
\begin{itemize}
\item In \cite{DeVMo} it was proved that {\em strictly positive vertex
    curvature} implies finiteness of a graph, thus proving a
  conjecture of Higuchi (which is a discrete analogue of Bonnet-Myers
  Theorem in Riemannian geometry). This question was investigated
  before by Stone \cite{St}.
\item The {\em cut locus} $\Cut(v)$ of a vertex $v$ consists of all
  vertices $w \in \V$, at which $d_v := d(v,\cdot)$ attains a local
  maximum, i.e., we have $w \in \Cut(v)$ if $d_v(w') \le d_v(w)$ for
  all $w' \sim w$. If $\G$ is a plane tessellation with {\em
    non-positive corner curvature}, then $\G$ is {\em without cut locus},
  i.e., we have $\Cut(v) = \emptyset$ for all $v \in \V$. This fact
  can be considered as a combinatorial analogue of the Cartan-Hadamard
  Theorem (for a proof and more details see \cite[Theorem 1]{BP2}).
\end{itemize}

For a finite subset $W \subset V$, let $\vol(W) = \sum_{v \in W}
|v|$. We consider the following two types of Cheeger constants:
\begin{equation} \label{eq:cheeger}
   \al(\G)=\inf_{\scriptsize\begin{array}{c} W\subseteq \V, \\
    |W|<\infty \end{array}} \frac{|\dd_E W|}{|W|}\quad\mathrm{and}\quad
    \alg(\G)=\inf_{\scriptsize\begin{array}{c} W\subseteq \V, \\
     |W|<\infty \end{array}} \frac{|\dd_E W|}{\vol(W)},
\end{equation}
where $\dd_E W$ is the set of all edges $e \in \E$ connecting a vertex
in $W$ with a vertex in $\V \backslash W$. The quantity $\al(\G)$ is
called the \emph{physical Cheeger constant} and $\alg(\G)$ the
\emph{geometric Cheeger constant} of the graph $\G$. The attributes
{\em physical} and {\em geometric} are motivated by the fact that
these constants are closely linked to two types of Laplacians (see,
e.g., \cite{Ke,We}) and that the first type is used in the community
of Mathematical Physics whereas the second appears frequently in the
context of Spectral Geometry. Cheeger constants are invariants of the
{\em global asymptotic geometry}. They are important geometric
tools for spectral considerations (both in setting of graphs and of
Riemannian manifolds) and play a prominent role in the topic of
expanders and Ramanujan graphs (see \cite{HLW} for a very
recommendable survey on this topic).

Natural model spaces are the $(p,q)$-regular plane tessellations
$\G_{p,q}$: every vertex in $G_{p,q}$ has degree $p$ and every face
has degree $q$. (In the case $\frac{1}{p} + \frac{1}{q} <
\frac{1}{2}$, $\G_{p,q}$ can be realised as a regular tessellation of
the Poincar{\'e} disc model of the hyperbolic plane by translates of a
regular compact polygon.) The graphs $\G_{p,q}$ can be considered as
discrete counterparts of constant curvature space forms in Riemannian
Geometry. The Cheeger constants of these regular graphs are
explicitely given:

\begin{Thm*}[see \cite{HJL,HiShi,HiShi2}]
  Let $\frac{1}{p} + \frac{1}{q} \le \frac{1}{2}$. Then
  $$ \alg(\G_{p,q}) = \frac{p-2}{p} \sqrt{1 - \frac{4}{(p-2)(q-2)}}. $$
\end{Thm*}

Let us now leave the situation of regular tessellations. It is known
that the Cheeger constants of general negatively curved planar graphs are
strictly positive (see \cite{Do}, \cite{Hi} and \cite{Woe}). Moreover, for
infinite planar graphs $G$ with $|v| \ge p$ and $|f| \ge q$ for almost
all vertices and faces and $c := \frac{1}{2} - \frac{1}{p} -
\frac{1}{q} > 0$, the following estimate was shown in \cite{Mo}:
\begin{equation} \label{eq:mocheeg}
\alpha(G) \ge \frac{2pqc}{3q-8}.
\end{equation}
Next we introduce a bit of notation before we state our explicit lower
Cheeger constant estimates. The variables $p,q$ in this paper always
represent a pair of numbers $3 \le p,q \le \infty$ satisfying
$\frac{1}{p} + \frac{1}{q} \le \frac{1}{2}$ (note that we use
$1/\infty =0$). For such a pair $(p,q)$, let
\begin{equation} \label{eq:cpq}
C_{p,q} := \begin{cases} 1, & \text{if $q= \infty$}, \\
1 + \frac{2}{q-2}, & \text{if $q < \infty$ and $p = \infty$}, \\
(1 + \frac{2}{q-2}) (1 + \frac{2}{(p-2)(q-2)-2}), & \text{if
$p,q < \infty$}. \end{cases}
\end{equation}
Then we have

\begin{Thm}[Cheeger constant estimate] \label{thm:cheegest}
  Let $\G = (\V,\E,\F)$ be a locally tessellating planar graph such that
  $|v| \le p\ \forall v \in \V$ and $|f| \le q\ \forall f \in \F$. (Note that
  $p = \infty$ or $q = \infty$ means no condition on the vertex of face
  degrees.) Let $C_{p,q}$ be defined as in \eqref{eq:cpq}.
  \begin{itemize}
  \item[(a)] Assume that $C := \inf_{v \in V} -\kappa(v)$ is strictly
    positive. Then
    $$ \al(\G) \ge 2 C_{p,q}C. $$
  \item[(b)] Assume that $c := \inf_{v \in V} -\frac{1}{|v|}\ka(v)$ is strictly
    positive. Then
    $$ \alg(\G) \ge 2 C_{p,q} c. $$
  \end{itemize}
  The above estimates are sharp in the case of regular trees (in which
  case $q = \infty$).
\end{Thm}

The proof of this theorem is given in Section \ref{sec:proofthm1}.
Observe that the constant $C_{p,q} \ge 1$ in \eqref{eq:cpq} becomes
largest if the graph $\G$ in Theorem \ref{thm:cheegest} has both finite
upper vertex and face degrees. (A shorter expression for $C_{p,q}$ is
$\frac{q(p-2)}{(p-2)(q-2)-2}$, which we have to interpret in the right way
if $q=\infty$ or $p=\infty$.)

Let us study our estimate in the regular case $\G = \G_{p,q}$: In
this case our estimate yields
$$ (p-2) \left( 1 - \frac{2}{(p-2)(q-2)-2} \right) \le \al(\G_{p,q}). $$
On the other hand, a straightforward calculation leads to the following
upper inequality
$$ \al(\G_{p,q}) = (p-2) \sqrt{1 - \frac{4}{(p-2)(q-2)}} \le (p-2) \left( 1 -
  \frac{2}{(p-2)(q-2)-1} \right), $$
which shows that our lower bound is very close to the correct
value. Mohar's estimate \eqref{eq:mocheeg} in this situation coincides
with ours in the particular case $(p,q) = (\infty,3)$, and becomes
considerably weaker for $q \ge 4$ or $p < \infty$.

\begin{Remark}
  Any infinite connected graph $\G = (\V,\E)$ with $|v| \le p$ has
  physical Cheeger constant $\alpha(\G) \le p-2$. To see this, choose
  an infinite path $v_0, v_1, v_2, \dots$ and let $W_n :=
  \{v_0,v_1,\dots,v_n\}$. Then we have
  $$ \frac{|\dd_E W_n|}{|W_n|} \le \frac{2(p-1)+(n-1)(p-2)}{n+1}, $$
  which implies
  $$ \al(\G) \le \lim_{n \to \infty} \frac{|\dd_E W_n|}{|W_n|} = p-2. $$
  The same arguments show $\al(\T_p) = p-2$ and $\alg(\T_p) =
  \frac{p-2}{2}$, where $\T_p$ denotes the $p$-regular infinite tree.
\end{Remark}

\smallskip

Next, we turn to another {\em global asymptotic invariants} related to
the growth of an infinite graph $\G = (\V,\E)$. For a fixed center
$v_0 \in \V$, let $S_n = S_n(v_0) = \{ v \in \V \mid d(v_0,v) = n \}$
be the spheres of radius $n$ and $\sigma_n = | S_n |$. The {\em growth
  series for $(\G,v_0)$} is the formal power series $f_{\G,v_0}(z) =
\sum_{n=0}^\infty \sigma_n z^n$ and the {\em exponential growth}
$\mu(\G,v_0)$ is given by
$$ \mu(\G,v_0) = \limsup_{n \to \infty} \frac{\log \sigma_n}{n}. $$
By Cauchy-Hadamard, the growth series represents a well-defined
function in the open complex ball of radius $e^{-\mu(\G,v_0)}$. In
many cases the exponential growth does not depend on the choice of
$v_0$. If this is the case, we simply write $\mu(\G)$.

Of particular importance in the study of the growth series
$f_{\G,v_0}$ are recursion formulas for the sequence $\sigma_n$. In
this paper, we consider the case of {\em $q$-face regular plane
  tessellations} $\G = (\V,\E,\F)$ (i.e., $|f|=q$ for all $f \in
\F$). Before stating our result in terms of average curvatures over
spheres $\overline{\kappa}(S_n) = \frac{\kappa(S_n)}{\sigma_n}$ we
need, again, some notation: For $3 \le q < \infty$ let $N = \frac{q-2}{2}$ if
$q$ is even and $N = q-2$ if $q$ is odd, and
\begin{equation} \label{eq:bl}
  b_l = \begin{cases} \frac{4}{q-2} & \text{if $q$ is even,}\\
    \frac{4}{q-2} & \text{if $q$ is odd and $l \neq \frac{N-1}{2}$,}\\
    \frac{4}{q-2}-2 & \text{if $q$ is odd and $l =
    \frac{N-1}{2}$,} \end{cases}
\end{equation}
for $0 \le l \le N-1$.

\begin{Thm}[Growth recursion formulas] \label{thm:facereggrow}
  Let $\G=(\V,\E,\F)$ be a $q$-face regular plane tessellation without
  cut locus, $S_n = S_n(v_0)$ for some $v_0 \in V$ and $\sigma_n =
  |S_n|$. Let $N$ and $b_l$ be defined as above (see
  \eqref{eq:bl}). Moreover, let $\curv_n =
  \frac{2q}{q-2}\overline{\kappa}(S_n)$. Then we have the following
  $(N+1)$-step recursion formulas for $n \ge 1$:
  \begin{equation} \label{eq:sigman}
  \sigma_{n+1} = \begin{cases} \sigma_1 + \sum_{l=0}^{n-1} (b_l - \curv_{n-l})
    \sigma_{n-l} & \text{if $n < N$,} \\[.3cm]
    \sum_{l=0}^{N-1} (b_l - \curv_{N-l}) \sigma_{N-l} & \text{if $n =
      N$,} \\[.3cm]
    - \sigma_{n-N} + \sum_{l=0}^{N-1} (b_l - \curv_{n-l}) \sigma_{n-l} &
    \text{if $n > N$.} \end{cases}
  \end{equation}
\end{Thm}

A proof of this theorem is given in Section \ref{sec:proofthm23}. Note
that the constants $\curv_k$ are zero for the regular flat tessellations
$\G_{3,6}, \G_{4,4}$ and $\G_{6,3}$. The constants $\curv_k$ in
\eqref{eq:sigman} can, therefore, be considered as curvature
correction terms for general non-flat tessellations.

\smallskip

In the special case of $(p,q)$-regular graphs $\G = \G_{p,q}$, the
terms $b_l - \curv_k$ all coincide with the constant $p-2$ except in the
case if $q$ is odd and $l = \frac{N-1}{2}$, when we have $b_l - \curv_k
= p-4$. In this case, Theorem \ref{thm:facereggrow} is equivalent to
the fact that $h_{p,q} f_{G,v_0} = g_{p,q}$ with
\begin{eqnarray*}
  h_{p,q} &=&  1 + 2z + \dots + 2z^N + z^{N+1},\\
  g_{p,q} &=& 1 - (p-2) z - \dots - (p-2) z^N + z^{N+1},
\end{eqnarray*}
if $q$ is even, and
\begin{eqnarray*}
  h_{p,q} &=& 1 + 2z + \dots + 2z^{\frac{N-1}{2}} + 4z^{\frac{N+1}{2}} +
  2z^{\frac{N+3}{2}} + \dots + 2z^N + z^{N+1}, \\
  g_{p,q} &=& 1 - (p-2)z - \dots - (p-4)z^{\frac{N+1}{2}} - \dots - (p-2)z^N
  + z^{N+1},
\end{eqnarray*}
if $q$ is odd. This agrees with results of Cannon and Wagreich
\cite{CaWa} and Floyd and Plotnick \cite[\S 3]{FP} that
the growth function $f_{\G,v_0}$ is the rational function
$h_{p,q}/g_{p,q}$. Moreover, it was shown in \cite{CaWa} and \cite{BaCS}
that the denominator polynomial $g_{p,q}$ for $\frac{1}{p} + \frac{1}{q} <
\frac{1}{2}$ is a reciprocal {\em Salem polynomial}, i.e., its roots lie on the
complex unit circle except for two positive reciprocal real zeros
$\frac{1}{x_{p,q}} < 1 < x_{p,q} < p-1$. This implies that the exponential growth
coincides with $\log x_{p,q}$, i.e.,
\begin{equation} \label{eq:muG}
\mu(\G_{p,q}) = \log x_{p,q} < \log (p-1) = \mu(\T_p).
\end{equation}
(An even more precise decription of the growth of the sequence
$\sigma_n$ is given in \cite[Cor. 3]{BaCS}.) Of course, it is
desirable to know more about the explicit value of $\mu(\G_{p,q}) =
\log x_{p,q}$. Since $g_{p,q}$ is divisible by $z^2 - (p-
\frac{4}{q-2}) z + 1$ in the case $q=3,4,6$, we have

\begin{Prop}
  Let $q \in \{3,4,6\}$. Then
  $$ \mu(\G_{p,q}) = \log \left( \frac{p}{2} - \frac{2}{q-2} + \sqrt{\left(
        \frac{p}{2} - \frac{2}{q-2} \right)^2 -1} \right). $$
\end{Prop}

In most of the cases the polynomial $g_{p,q}$ is essentially
irreducible (expect for some small well known factors; see
\cite[Thm. 1]{BaCS}) and there is no hope to have an explicit
expression for its largest zero $x_{p,q} > 1$. A direct consequence of
the isoperimetric inequality in \cite[Cor. 5.2]{BP1} is the following
lower estimate of $\log x_{p,q}$:

\begin{Prop} \label{prop:lpqest}
  Let $\G_{p,q}$ be non-positively curved, i.e.,
  $$ C = -\kappa(v) = p \left( \frac{1}{p} + \frac{1}{q} - \frac{1}{2}
  \right) \ge 0 \quad \forall \ v \in \V, $$
  then we have
  \begin{equation} \label{eq:lpqest}
  \mu(\G_{p,q}) = \log x_{p,q} \ge \log \left( 1 + \frac{2q}{q-1} C \right).
  \end{equation}
\end{Prop}

Note that \eqref{eq:lpqest} implies $\lim_{q \to \infty} \mu(\G_{p,q})
= \mu(\T_p) = \log p-1$ for all $p \ge 3$.

\begin{Remark}
  The {\em Mahler measure} $M(g)$ of a monic polynomial $g \in \Z[z]$
  with integer coefficients is given by the product $\prod |z_i|$,
  where $z_i \in \C$ are the roots of $g$ of modulus $\ge 1$. {\em
    Lehmer's conjecture} states that for every such $g$ with $M(g) >
  1$ we have
  $$ M(g) \ge M(1-z+z^3-z^4+z^5-z^6+z^7-z^9+z^{10}) \approx 1.1762\dots. $$
  Thus \eqref{eq:lpqest} yields an explicit lower estimate for the
  Mahler measure of the polynomials $g_{p,q}(z)$.
\end{Remark}

Let us now return to general $q$-face regular tessellations. We conclude from
Theorem \ref{thm:facereggrow}:

\begin{Thm}[Curvature/Growth comparison] \label{thm:growest}
  Let $\G = (\V,\E,\F)$ and $\widetilde \G = (\widetilde \V,\widetilde
  \E, \widetilde \F)$ be two $q$-face regular plane tessellations with
  non-positive vertex curvature, $S_n \subset \V$ and $\widetilde S_n
  \subset \V$ be spheres with respect to the centres $v_0 \in \V$ and
  $\widetilde v_0 \in \widetilde \V$, respectively, and $\sigma_n =
  |S_n|$ and $\widetilde \sigma_n = |\widetilde S_n|$. Assume that the
  spherical average curvatures satisfy
  $$ \overline{\kappa}(\widetilde S_n) \le \overline{\kappa}( S_n) \le 0 \quad
  \forall \ n \ge 0. $$
  Then the difference sequence $\widetilde \sigma_n - \sigma_n \ge 0$
  is monotone non-decreasing and, in particular, we have $\mu(\widetilde \G) \ge
  \mu(\G)$.
\end{Thm}

A proof of this theorem is given in Section \ref{sec:proofthm23}. (In
fact, the proof shows that the vertex curvature conditions in Theorem
\ref{thm:growest} can be slightly relaxed: It suffices that
$\G$ has non-positive vertex curvature and that both graphs
$\G$ and $\widetilde \G$ are without cut-loci.) This theorem can be
considered as a refined discrete counterpart of the
Bishop-G\"unther-Gromov Comparison Theorem for Riemannian manifolds
(see, e.g., \cite[Theorem 3.101]{GaHuLa}). The latter compares volumes
of balls in Riemannian manifolds against constant curvature space
forms; our discrete counterpart deals with spheres (the result for
balls is obtained by adding over spheres) and is more flexible as it
allows to use more general comparison spaces.

If we drop the face regularity condition, it is not difficult to
derive the following simple tree comparison result. A proof of this
result is given in Section~\ref{sec:proofthm4}:

\begin{Thm}[Tree comparison] \label{thm:treecomp} Let $\G =
  (\V,\E,\F)$ be a locally tessellating planar graph without cut locus
  satisfying $|v| \le p$ for all $v \in \V$, for some $p \ge 3$. Then
  $$ \mu(\G) \le \mu(\T_p) = \log(p-1). $$
\end{Thm}

Note that another more involved tree comparison result was obtained by
Higuchi \cite{Hi2} for (not necessarily planar) infinite
vertex-regular graphs with each vertex contained a cycle of uniformly
bounded length.

\smallskip

Let us finally discuss some spectral applications (see \cite{MW,Woe2}
for classical surveys). The (geometric) Laplacian
$\Lp:\ell^2(V,m)\lra\ell^2(V,m)$, where $m(v)=|v|$ for $v\in V$ is
given by
$$(\Lp\ph)(v)=\frac{1}{|v|}\sum_{w\sim v}\ph(v)-\ph(w),
\qquad(\ph\in\ell^2(V,m),v\in V).$$
The relation between the bottom $\lm_0(\G)$ of the spectrum and the bottom $\lm_0^\ess(\G)$ of the essential spectrum of
$\Lp$ and the Cheeger constant and the exponential growth in the discrete case was presented first by Dodziuk/Kendall \cite{DKe} and Dodziuk/Karp \cite{DKa}. The best estimates are
due to K. Fujiwara (see \cite{Fu1} and \cite{Fu2}):
\begin{equation} \label{eq:fuji}
1-\sqrt{1-\widetilde\al^2(\G)} \le \lambda_0(\G) \le
\lambda_0^\ess(\G) \le 1 - \frac{2 e^{\mu(\G)/2}}{1+e^{\mu(\G)}},
\end{equation}
which are sharp in the case of regular trees.

An immediate consequence of Theorem \ref{thm:cheegest} and
\eqref{eq:fuji} is the following combinatorial analogue of McKean's
Theorem (see \cite{McK} for the result in the smooth setting):

\begin{Cor}
  Let $\G = (\V,\E,\F)$ be a locally tessellating planar graph
  satisfying the vertex and face degree bounds in Theorem
  \ref{thm:cheegest}. Moreover, assume that $c := \inf_{v \in V} -
  \frac{1}{|v|}\kappa(v) > 0$. Then
  $$ 1-\sqrt{1-\left( 2 C_{p,q} c \right)^2} \le \lambda_0(\G), $$
  where $C_{p,q}$ is defined in \eqref{eq:cpq}. This estimate is sharp
  in the case of regular trees.
\end{Cor}

Similarly, Theorem \ref{thm:growest}, \eqref{eq:muG} and
\eqref{eq:fuji} directly imply

\begin{Cor}
  Let $p,q \ge 3$ and $\frac{1}{p} + \frac{1}{q} \le \frac{1}{2}$.
  Let $\G$ be a $q$-face regular tessellation without cut locus and
  satisfying $|v| \le p$ for all vertices. Then
  $$ \lambda_0^\ess(\G) \le 1 - \frac{2 \sqrt{x_{p,q}}}
  {1+x_{p,q}}. $$
  This estimate is sharp in the case of regular trees.
\end{Cor}

Let us finish this introduction with some general references. It was
shown in \cite{KLPS} that non-positive corner curvature implies
non-existence of finitely supported eigenfunctions of all elliptic
operator on planar graphs.  Lower estimates for the bottom of the
essential spectrum in terms of {\em Cheeger constants at infinity} or
branching rates of general non-planar graphs can be found in
\cite{Fu2} and \cite{Ura} for the geometric Laplacian and in \cite{Ke}
and \cite{Woj,Woj2} for the physical Laplacian
$\Delta:D(\Delta)\subseteq\ell^2(V)\lra\ell^2(V)$ given by
$(\Delta\ph)(v)=\sum_{w\sim v}\ph(v)-\ph(w)$, $\ph\in D(\Delta),v\in
V$. These results show, in particular, absence of the essential
spectrum for graphs with curvature converging to infinity outside
increasing compact sets, a phenomenon which was first proved in the
smooth context of manifolds by \cite{DL}.

\medskip

{\bf Acknowledgement:} Matthias Keller likes to thank Daniel Lenz who
encouraged him to study the connection between curvature and spectral
theory. Matthias Keller was supported during this work by sdw.
Norbert Peyerimhoff is grateful for the financial support of the
TU Chemnitz. Both authors like to thank Ruth
Kellerhals and Victor Abrashkin for very useful discussions.

\section{Proof of Theorem \ref{thm:cheegest}}
\label{sec:proofthm1}

Let us first give precise definitions of some notions used in the
introduction. Let $\G=(\V,\E)$ be a planar graph embedded in
$\R^2$. The faces $f$ of $\G$ are the closures of the connected
components in $\R^2\setminus \bigcup_{e\in E} e$.

We further assume that $\G$ has no loops, no multiple edges and no
vertices of degree one (leaves). Moreover, we assume that
every vertex has finite degree and that every bounded open set in
$\R^2$ meets only finitely many faces of $\G$. We call a planar graph
with these properties {\em simple}. We call a sequence of edges
$e_1,\dots,e_n$ a {\em walk of length $n$} if there is a corresponding
sequence of vertices $v_1,\dots,v_{n+1}$ such that $e_i =
v_iv_{i+1}$. A walk is called a {\em path} if
there is no repetition in the corresponding sequence of vertices
$v_1,\dots,v_n$. A (finite or infinite) path with associated vertex
sequence $\dots v_i v_{i+1} v_{i+2} \dots$ is called a {\em geodesic},
if we have $d(v_i,v_j) = |i-j|$ for all pairs of vertices in the path.
The {\em boundary of a face $f$} is the subgraph
$\partial f = (\V \cap f,\E \cap f)$. We define the degree $|f|$ of a
face $f\in F$ to be the length of the shortest closed walk in the
subgraph $\partial f$ meeting all its vertices. If there is no such
finite walk, we set $|f|=\infty$. Now we present the conditions which have
to be satisfied that a planar graph is locally tessellating:

\begin{Definition} \label{def:loctessgraph}
A simple planar graph $\G$ is called a {\em locally tessellating
  planar graph} if the following conditions are satisfied:
\begin{itemize}
\item [(i)] Any edge is contained in precisely two different faces.
\item [(ii)] Any two faces are either disjoint or have precisely a
  vertex or a path of edges in common. In the case that the length of
  the path is greater than one, then both faces are unbounded.
\item [(iii)] Any face is homeomorphic to the closure of an open disc
  ${\mathbb D} \subset \R^2$, to $\R^2\setminus {\mathbb D}$ or to the
  upper half plane $\R\times\R_+\subset \R^2$ and its boundary is a
  path.
\end{itemize}
\end{Definition}

Note that these properties force the graph $\G$ to be connected.
Examples are tessellations $\R^2$ introduced in \cite{BP1,BP2}, trees
in $\R^2$, and particular finite tessellations on the sphere mapped to
$\R^2$ via stereographic projection.

\medskip

Now we turn to the proof of Theorem \ref{thm:cheegest}. The heart of
the proof is Proposition \ref{prop:Harm} below. An earlier version of
this proposition in the dual setting is Proposition 2.1 of \cite{BP1}.
We start with a few more preliminary considerations. Let
$\G=(\V,\E,\F)$ be a locally tessellating planar graph. For a finite
set $W\subseteq \V$ let $\G_W=(W,\E_W,\F_W)$ be the subgraph of $\G$
induced by $W$, where $\E_W$ are the edges in $\E$ with both end
points in $W$ and $\F_W$ are the faces induced by the graph
$(W,\E_W)$. Euler's formula states for a finite and connected subgraph
$\G_W$ (observe that $\F_W$ contains also the unbounded face):
\begin{equation}\label{eq:Euler}
  |W|-|\E_W|+|\F_W|=2.
\end{equation}
Recall that $\dd_E W$ is the set of edges connecting a vertex in $W$
with one in $V \backslash W$. By $\dd_F W$, we denote the set of faces
in $F$ which contain an edge of $\dd_E W$. Moreover, we define the
\emph{inner degree} of a face $f\in\dd_F W$ by
$$|f|^i_W=|f\cap W|.$$

We will need two useful formulas which hold for
arbitrary finite and connected subgraphs $\G_W=(W,\E_W,\F_W)$. The first
formula is easy to see and reads as
\begin{equation}\label{eq:E_W}
  \sum_{v\in W}|v|=2|\E_W|+|\dd_E W|.
\end{equation}
Since $W$ is finite, the set $\F_W$ contains at least one face which is
not in $\F$, namely the unbounded face surrounding $\G_W$, but there can
be more. Define $c(W)$ as the number
\begin{equation} \label{eq:CW}
c(W)=|\F_W|-|\F_W\cap \F| \ge 1.
\end{equation}
Note that $|\F_W \cap \F|$
is the number of faces in $\F$ which are entirely enclosed by edges of
$\E_W$. Sorting the following sum over vertices according to faces
gives the second formula
\begin{eqnarray}
\sum_{ v\in W}\sum_{f\ni v}\frac{1}{|f|} &=& |\F_W \cap \F| +
\sum_{f\in\dd_F W}\frac{{|f|}^i_W}{|f|} \nonumber \\
&=& |\F_W| - c(W) + \sum_{f\in\dd_F W}\frac{{|f|}^i_W}{|f|}. \label{eq:F_W}
\end{eqnarray}

\begin{Prop}\label{prop:Harm}
  Let $\G=(\V,\E,\F)$ be a locally tessellating planar graph and $W
  \subset \V$ be a finite set of vertices such that the induced
  subgraph $\G_W$ is connected. Then we have
  $$\ka(W)=2 -c(W) -\frac{|\dd_E W|}{2}+\sum_{f\in\dd_F W}\frac{{|f|}^i_W}{|f|}.$$
\end{Prop}

\begin{proof}
  By the equations \eqref{eq:E_W}, \eqref{eq:F_W} and \eqref{eq:Euler}
  we conclude
  \begin{eqnarray*}
    \ka(W)&=&\sum_{ v\in W}\left(1-\frac{|v|}{2}+\sum_{f\ni
        v}\frac{1}{|f|}\right)\\
    &=&|W|-{|\E_W|}-\frac{|\dd_E W|}{2}+{|\F_W|}-c(W)+
    \sum_{f\in\dd_F
      W}\frac{{|f|}^i_W}{|f|}\\
    &=& 2 -c(W) -\frac{|\dd_E W|}{2}+\sum_{f\in\dd_F W}\frac{{|f|}^i_W}{|f|}.
  \end{eqnarray*}
\end{proof}

A finite set $W \subset V$ is called a {\em polygon}, if $G_W$ is
connected and if $c(W) = 1$. This notion becomes understandable
if one looks {\em at the dual setting}: Every vertex $v \in W$
corresponds to a face $f^*(v) \in \F^*$ in the dual planar graph $\G^*
= (\V^*,\E^*,\F^*)$, and $W \subset \V$ is a polygon if and only if
$\bigcup_{v \in W } \overline{f^*(v)} \subset \R^2$ is homeomorphic to
a closed disc (here $\overline{f}$ denotes the closure of the
geometric realization of the face $f$). For $v \in W$, let $|v|_W^e$
denote the number of edges in $\dd_E W$ adjacent to $v$. $|v|_w^e$ is called
the {\em external degree} of $v$ (w.r.t. $W$). Moreover, let
$\dd_V W$ be the set of vertices in $W$ with $|v|^e_W \ge 1$.

\begin{Prop} \label{prop:isop}
  Ler $\G=(\V,\E,\F)$ be a locally tessellating planar graph
  satisfying the vertex and face degree bounds in Theorem
  \ref{thm:cheegest} and $W \subset \V$ be a polygon with $|v|^e_W \le
  p-2$ for all $v \in \dd_V W$. Then we have
  \begin{equation} \label{eq:ddEWest}
  |\dd_E W| \ge 2 C_{p,q} (1-\kappa(W)).
  \end{equation}
  Moreover, under the assumption of (a) or (b) in Theorem \ref{thm:cheegest},
  we have
  $$ \frac{|\dd_E W|}{|W|} \ge 2 C_{p,q} C \quad \text{or} \quad
  \frac{|\dd_E W|}{\vol(W)} \ge 2 C_{p,q} c, \ \text{respectively.} $$
\end{Prop}

\begin{proof}
  Observe first that we have the inequality
  \begin{equation} \label{eq:fiW}
  \sum_{f\in\dd_F W} {{|f|}^i_W} \ge |\dd_V W| + |\dd_E W|.
  \end{equation}
  This can be seen as follows: Every face $f \in \dd_F W$ may have
  some edges and some isolated vertices in common with the induced
  graph $\G_W=(W,\E_W,\F_W)$. Since the vertices of $\dd_V W$ are
  connected in $\G_W$, there are at least $|\dd_V W|$ pairs $(f,e) \in
  \dd_F W \times E$ with $e \in \partial f \cap E_W$. These pairs
  contribute at least $2 |\dd_V W|$ to the left hand sum in
  \eqref{eq:fiW}. At every vertex $v \in \dd_V W$, there are
  $|v|_W^e-1$ faces of $\dd_F W$ which meet $G_W$ in the isolated
  vertex $v$. Adding over all these vertices $v \in \dd_V W$, we
  obtain the total contribution $|\dd_E W| - |\dd_V W|$ to the left
  hand sum in \eqref{eq:fiW}. One easily checks that there is no
  overlap of both contributions, leading to the above inequality.

  Using \eqref{eq:fiW}, $|\dd_V W| \ge \frac{1}{p-2}|\dd_E W|$, $|f|
  \le q$ for all $f \in \F$, and Proposition \ref{prop:Harm}, we
  obtain
  $$ |\dd_E W| \left( \frac{1}{2} - \frac{p-1}{q(p-2)} \right) \ge
  \frac{|\dd_E W|}{2} - \sum_{f\in\dd_F W} \frac{{|f|}^i_W}{|f|} = 1-
  \kappa(W), $$
  which yields \eqref{eq:ddEWest}. The second formula of the
  proposition follows from $-\kappa(W) \ge C |W|$ in case (a) and from
  $-\kappa(W) \ge c \vol(W)$ in case (b).
\end{proof}

Henceforth, let $\G = (\V,\E,\F)$ be a locally tessellating planar
graph as in Theorem \ref{thm:cheegest}. Recall that $C_{p,q} =
\frac{q(p-2)}{(p-2)(q-2)-2}$. The conditions $|v| \le p$ and $|f| \le
q$ for all $v \in \V$ and $f \in \F$ imply $2 C_{p,q} C \le p-2$ and
$2 C_{p,q} c \le \frac{p-2}{p}$.

\begin{Lemma} \label{lem:help1}
  Let $v \in V$ and $W = \{ v \}$. Then we have
  \begin{equation} \label{eq:alv}
  \frac{|\dd_E W|}{|W|} = |v| \ge 2 C_{p,q} C
  \end{equation}
  and
  $$
  \frac{|\dd_E W|}{\vol(W)} = 1 \ge \frac{p-2}{p}.
  $$
\end{Lemma}

\begin{proof}
  The only non trivial inequality is \eqref{eq:alv}. It follows
  straightforwardly from $\kappa(v) \le - C$ that
  $$ C \le \frac{q-2}{2q}|v| - 1. $$
  This implies that
  \begin{eqnarray*}
    2 C_{p,q} C &\le& \frac{(p-2)(q-2)}{(p-2)(q-2)-2} |v| - \frac{2q(p-2)}
    {(p-2)(q-2)-2} \\
    &\le& |v| - \frac{2}{(p-2)(q-2)-2} (q(p-2)-p).
  \end{eqnarray*}
  The lemma follows now from the fact that $q(p-2)-p \ge 0$ for $p,q \ge 3$.
\end{proof}

\begin{Lemma} \label{lem:help2}
  Assume that there is a finite set $W \subset V$ such that
  \begin{equation} \label{eq:isopcond}
    \frac{|\dd_E W|}{|W|} < 2 C_{p,q} C \le p-2 \quad \text{or} \quad
    \frac{|\dd_E W|}{\vol(W)} < 2 C_{p,q} c \le \frac{p-2}{p}, \
    \text{respectively.}
  \end{equation}
  Then there exists a {\em polygon} $W' \subset V$ with $|v|_{W'}^e
  \le p-2$ for all $v \in \dd_V W'$, such that
  $$ \frac{|\dd_E W'|}{|W'|} \le \frac{|\dd_E W|}{|W|} \quad \text{or} \quad
  \frac{|\dd_E W'|}{\vol(W')} \le \frac{|\dd_E W|}{\vol(W)}, \
  \text{respectively.} $$
\end{Lemma}

\begin{proof}
  Observe first that we can always find a non-empty subset $W_0
  \subset W$ such that $G_{W_0}$ is a connected component of $G_W$ and
  that $|\dd_E W_0|/|W_0| \le |\dd_E W|/|W|$ or $|\dd_E W_0|/\vol(W_0)
  \le |\dd_E W|/\vol(W)$, respectively. Note that $G_{W_0}$ has only
  one unbounded face. By adding all vertices of $V$ contained in the
  union of all bounded faces of $G_{w_0}$, we obtain a polygon with
  even smaller isoperimetric constants. Let us denote this non-empty
  polygon, again, by $W$. By Lemma \ref{lem:help1}, $W$ must have at
  least two vertices. By connectedness of $G_W$ and $|W| \ge 2$, we
  have $|v|_W^e \le p-1$ for all $v \in W$. Assume there is a vertex
  $v \in \dd_VW$ with $|v|_W^e = p-1$.  Let $W' := W \backslash \{ v
  \}$. Then one easily checks that the condition \eqref{eq:isopcond}
  implies
  $$ \frac{|\dd_E W'|}{|W'|} = \frac{|\dd_E W|+2-p}{|W|-1} < \frac{|\dd_E W|}
  {|W|}$$
  or
  $$ \frac{|\dd_E W'|}{\vol(W')} = \frac{|\dd_E W|+2-p}{\vol(W)-p} <
  \frac{|\dd_E W|}{\vol(W)}, $$ respectively. Repeating this
  elimination of vertices with external degree $p-1$, we end up with a
  polygon $W'$ satisfying $|v|_{W'}^e \le p-2$ for all $v \in W'$ or
  with $W'$ equal to a single vertex. But the latter case is a
  contradiction to Lemma \ref{lem:help1}.
\end{proof}

\begin{proof}[Proof of Theorem \ref{thm:cheegest}] Since $2 C_{p,q} C \le p-2$ or $2 C_{p,q} c \le
\frac{p-2}{p}$, we only have to consider the cases when $\al(\G) <
p-2$ or $\alg(\G) < \frac{p-2}{p}$, since otherwise there is nothing
to prove. Lemma \ref{lem:help2} states that, if there is a finite $W
\subset V$ with
$$ \frac{|\dd_E W|}{|W|} < 2 C_{p,q} C \quad \text{or} \quad
\frac{|\dd_E W|}{\vol(W)} < 2 C_{p,q} c, \ \text{respectively,} $$
then there is a polygon $W'$ with $|v|_{W'}^e \le p-2$ for all $v \in
\dd_VW'$ satisfying the same inequality. But this contradicts to
Proposition \ref{prop:isop}, finishing the proof of Theorem
\ref{thm:cheegest}. \end{proof}

\section{Proofs of Theorems \ref{thm:facereggrow} and
  \ref{thm:growest}}
\label{sec:proofthm23}

Let $\G=(\V,\E,\F)$ be a $q$-face regular plane tessellation without
cut locus, $v_0 \in V$, $S_n = S_n(v_0)$ and $\sigma_n =
|S_n|$. Recall that we have $N = \frac{q-2}{2}$ if $q$ is even and $N
= q-2$ if $q$ is odd.  The recursion formulas \eqref{eq:sigman} in
Theorem \ref{thm:facereggrow} for $n \le N$ and $n > N$, respectively,
require separate proofs. However, both proofs are based on the
following results from \cite[Section 6]{BP1} (note that these results
are presented there in the dual setting of vertex-regular graphs):
Proposition 6.3 in \cite{BP1} states for $n \ge 1$ that
\begin{equation}\label{eq:ka}
  \kappa(B_n) = 1 - \frac{q-2}{2q}(\sigma_{n+1}-\sigma_n) +
  \sum_{j=2}^{q-2}\frac{q-2j}{2q} c_n^j,
\end{equation}
where $B_n = \{ v \in V \mid d(v_0,v) \le n \}$ denotes the ball and
$$c_n^j=\left| \{f \in \F \mid | f \cap ( \V \backslash B_n) | = j \}
\right|$$
for $1 \leq j \leq q-1$. Moreover from Lemma 6.2 in
\cite{BP1} we have the following recurrence relations for $c_n^j$, $n \geq 1$,
which arise very naturally from the geometric context
\begin{itemize}
  \item[(i)] $c_{n}^{l}=c_{n-1}^{l+2}$, for $1\leq l\leq q-3$,
  \item[(ii)] $c_{n}^{q-2}=c_{n-1}^2$,
  \item[(iii)] $c_{n}^{q-1}=c_{n}^1+\sigma_{n+1}-\sigma_{n}=
   c_{n-1}^3+\sigma_{n+1}-\sigma_{n}$.
\end{itemize}

We first aim at the proof of \eqref{eq:sigman} for $n \le N$. Let
$\tau: \{0,1,\dots,N\} \to \{1,2,\dots,q-1\}$ be defined as
$$
\tau(k) = \begin{cases} q-1-2k & \text{if $q$ is even and
    $0 \le k \le N$,}\\
  q-1-2k & \text{if $q$ is odd and $0 \le k \le \frac{N-1}{2}$,}\\
  2q-3-2k & \text{if $q$ is odd and $\frac{N+1}{2} \le k \le N$.}
\end{cases}
$$
Note that $\tau$ is defined precisely in such a way that we have
\begin{equation} \label{eq:cntau}
c_{n+1}^{\tau(k+1)} = c_n^{\tau(k)} \quad \text{for $0 \le k \le N-1$ and
$n \ge 0$,}
\end{equation}
by the recurrence relations (i) and (ii).

\begin{Lemma} \label{lem:cntau}
  Let $1 \le n \le N$. Then we have
  $$ c_n^{\tau(l)} = \begin{cases} \sigma_{n+1-l} - \sigma_{n-l} & \text{for $1
      \le l \le n-1$,}\\
    \sigma_1 & \text{for $l=n$},\\
    0 & \text{for $n+1 \le l \le N$.} \end{cases}
  $$
  Moreover, in the case $q$ even, we have $c_n^{2l} = 0$ for $1 \le l
  \le N$.
\end{Lemma}

\begin{proof}
  One easily sees that $c_0^j =0$ for $1 \le j \le q-2$ and
  $c_0^{q-1}=\sigma_1$. The recurrence relations (i) and (ii) imply that
  $c_k^1 = 0$ for $0 \le k \le N-1$ and $c_N^1 = \sigma_1$. Using
  (iii), we obtain $c_k^{q-1} = \sigma_{k+1}-\sigma_k$ for $1 \le k
  \le N-1$.  The value of $c_n^{\tau(l)}$ can now be deduced from
  these results by repeatedly applying \eqref{eq:cntau} in each of the
  cases $1 \le l \le n-1$, $l=n$ and $n+1 \le l \le N$.
\end{proof}

\begin{Lemma} \label{lem:recrel}
  Let $N \ge 2$, $1 \le n \le N$ and $b_l$ be defined as in
  \eqref{eq:bl}. Then we have
  $$
    \frac{q-6}{q-2} \sigma_n + \sum_{j=2}^{q-2} \frac{q-2j}{q-2} c_n^j
    = \begin{cases} \sum_{l=1}^{n-1} b_l \sigma_{n-l}
      & \text{if $n \le N-1$,} \\
      \frac{6-q}{q-2} \sigma_1 +\sum_{l=1}^{N-2} b_l \sigma_{N-l}
      & \text{if $n=N$.} \end{cases}
  $$
\end{Lemma}

\begin{proof}
  First observe that, since $2 \le \tau(l) \le q-2$ for
  $1 \le l \le N-1$ and $c_n^{\tau(l)}=0$ for $n+1 \le l \le N$ by Lemma \ref{lem:cntau}
  $$ \sum_{j=2}^{q-2} (q-2j) c_n^j
  = \begin{cases} \sum_{l=1}^n (q-2\tau(l)) c_n^{\tau(l)} & \text{if
      $n \le N-1$,} \\ \sum_{l=1}^{N-1} (q-2\tau(l)) c_n^{\tau(l)} &
    \text{if $n=N$,} \end{cases}
  $$
  (Note for the case $n=N$: we have $2 \le \tau(l) \le q-2$ only for
  $1 \le l \le N-1$ and $\tau(N) = 1$. This makes it necessary to
  treat this case separately.) The proof follows now straightforwardly
  with the help of Lemma \ref{lem:cntau} and the equation $(q-2) b_l =
  2(\tau(l) - \tau(l+1))$.
\end{proof}

\begin{proof}[Proof of Theorem \ref{thm:facereggrow}]
We rewrite equation \eqref{eq:ka} as follows:
$$ \sigma_{n+1} - \sigma_n = \frac{2q}{q-2} - \frac{2q}{q-2} \sum_{l=0}^n
\kappa(S_l) + \sum_{j=2}^{q-2} \frac{q-2j}{q-2} c_n^j. $$
Using $\kappa(S_0) = 1 - \frac{q-2}{2q}\sigma_1$ and
$\frac{2q}{q-2}\kappa(S_l) = \curv_l \sigma_l$ (for the definition of
$\curv_l$ see Theorem \ref{thm:facereggrow}) we obtain
\begin{equation} \label{eq:sigmanp1}
\sigma_{n+1} = \left( \sigma_1 - \sum_{l=1}^n \curv_l \sigma_l \right) +
\sigma_n + \sum_{j=2}^{q-2} \frac{q-2j}{q-2} c_n^j.
\end{equation}

The recursion formulas in Theorem \ref{thm:facereggrow} for $N \ge 2$
and $n \le N$ follow now directly from \eqref{eq:sigmanp1} and Lemma
\ref{lem:recrel}. The case $n=N=1$ has to be treated separately: In
this case we have $\sum_{j=2}^{q-2} (q-2j) c_n^j = 0$ and
\eqref{eq:sigmanp1} simplifies to $\sigma_2 = (2-\curv_1)\sigma_1$. The result
follows now from the fact that $b_0 = 2$.

It remains to prove the recursion formula \eqref{eq:sigman} for $n >
N$. We first consider the case $N \ge 2$. Repeated application of the
recurrence relations (i)-(iii) yields
\begin{equation} \label{eq:relN}
\sum_{j=2}^{q-2} \frac{q-2j}{q-2} c_n^j = \left( \sum_{l=0}^{N-1}
b_l \sigma_{n-l}\right) - (\sigma_n + \sigma_{n-(N-1)}) + \sum_{j=2}^{q-2}
\frac{q-2j}{q-2} c_{n-N}^j.
\end{equation}
Since $\frac{2q}{q-2} (\kappa(B_n) - \kappa(B_{n-N})) =
\sum_{l=0}^{N-1} \curv_{n-l} \sigma_{n-l}$, we obtain with \eqref{eq:ka}
and \eqref{eq:relN}
\begin{multline*}
  \sum_{l=0}^{N-1} \curv_{n-l} \sigma_{n-l} = \\ = -(\sigma_{n+1}
  -\sigma_n) + (\sigma_{n-(N-1)} - \sigma_{n-N}) + \left(
    \sum_{l=0}^{N-1} b_l \sigma_{n-l} \right) - (\sigma_n +
  \sigma_{n-(N-1)}),
\end{multline*}
which immediately yields the recursion formula \eqref{eq:sigman} for
$n > N \ge 2$. The case $N=1$ is particularly easy and left to the
reader.  This finishes the proof of Theorem
\ref{thm:facereggrow}.
\end{proof}

\medskip

Now we turn to the proof of Theorem \ref{thm:growest}. Note first that non-positive
vertex curvature implies non-positive corner curvature in the case of
face-regular graphs. By \cite[Theorem 1]{BP2}, both graphs $G,
\widetilde G$ are therefore without cut-loci and we can apply the
recursion formulas in Theorem \ref{thm:facereggrow}.

\begin{Lemma} \label{lem:blmimuk}
  We have the following estimates for $0 \le l \le N-1$ and $k \ge 1$:
  \begin{align*}
    b_l -  \curv_k &\ge 2, && \text{if $l = \frac{N-1}{2} (= 0)$ and
      $q \in \{3,4\}$,}\\
    b_l -  \curv_k &\ge 1, && \text{if $l \neq \frac{N-1}{2}$ or
      $q$ even,}\\
    b_l -  \curv_k &\ge 0, && \text{if $l = \frac{N-1}{2} (= 1)$ and
      $q=5$,}\\
    b_l -  \curv_k &\ge -1, && \text{if $l = \frac{N-1}{2}$ and
      $q \ge 7$ odd.}
  \end{align*}
\end{Lemma}

\begin{proof}
  The case ``$l \neq \frac{N-1}{2}$ or $q$ even'' follows from $b_l =
  \frac{4}{q-2}$ and $\kappa(v) = 1 - \frac{q-2}{2q}|v| \le 1 - 3
  \frac{q-2}{2q}$.

  Now assume that $l = \frac{N-1}{2}$. Since $b_l \ge \frac{4}{q-2} -
  2$, the previous considerations lead to $b_l -  \curv_k \ge
  -1$. If $q=3$ or $q=4$, then $b_l = 2$ and, consequently, $b_l -
   \curv_k \ge b_l = 2$. Finally, if $q=5$, then $\kappa(v)
  \le 0$ implies that $|v| \ge 4$ and thus $\kappa(v) \le 1 - 4
  \frac{q-2}{2q}$. Using this fact leads directly to $b_l -
  \curv_k \ge 0$.
\end{proof}

From the above lemma we deduce the following facts:

\begin{Lemma} \label{lem:4facts}
  We have
  \begin{itemize}
  \item[(a)] $b_0 -  \curv_k \ge 1$.
  \item[(b)] $b_{N-1} - \curv_k \ge 1$ if $N \ge 2$.
  \item[(c)] $b_0 - \curv_k \ge 2$ if $q=3$ or $q=4$.
  \item[(d)] Let $n,N \ge 1$, $1 \le k \le \min \{n,N\}$, and assume
    that $\gamma_i \ge 0$ are mononote non-decreasing for $n-k+1 \le i
    \le n-1$. Then
    $$ \sum_{l=1}^{k-1} (b_l -  \curv_{n-l}) \gamma_{n-l} \ge 0. $$
  \end{itemize}
\end{Lemma}

\begin{proof}
  (a), (b) and (c) are trivial consequences of Lemma
  \ref{lem:blmimuk}. (d) follows immediately from Lemma
  \ref{lem:blmimuk}, unless we have $k-1 \ge \frac{N-1}{2}$ and $q \ge
  7$ odd. But in this case we have $\frac{N-1}{2} \ge 2$ and
  \begin{eqnarray*}
    \sum_{l=1}^{k-1} (b_l -  \curv_{n-l}) \gamma_{n-l} &\ge&
    (b_1 -  \curv_{n-1}) \gamma_{n-1} + (b_{\frac{N-1}{2}} -
     \curv_{n-\frac{N-1}{2}}) \gamma_{n-\frac{N-1}{2}}\\
    &\ge& \gamma_{n-1} - \gamma_{n-\frac{N-1}{2}} \ge 0,
    \end{eqnarray*}
    by the monotonicity of $\gamma_i$.
\end{proof}

\begin{proof}[Proof of Theorem \ref{thm:growest}] The condition
$\overline{\kappa}(\widetilde S_n) \le \overline{\kappa}(S_n) \le 0$
implies that $\widetilde \curv_n \le \curv_n$ for $n \ge 0$, where $\widetilde \curv_n=
  \frac{2q}{q-2}\overline{\kappa}(\widetilde S_n)$. Note that
$\sigma_0 = \widetilde \sigma_0 = 1$. Since $\sigma_1 = \frac{2q}{q-2}
- \curv_0$, $\widetilde\sigma_1 = \frac{2q}{q-2}
- \widetilde\curv_0$ we conclude that
$$ \widetilde \sigma_1 - \sigma_1 =\curv_0 -  \widetilde \curv_0 \ge 0 =
\widetilde \sigma_0 - \sigma_0. $$

The proof is based on induction over $n$: Assume that $n \ge 1$ and
that $\gamma_k :=  \widetilde \sigma_k -\sigma_k$ is non-negative and
monotone non-decreasing for $0 \le k \le n$. We aim to show that
$\gamma_{n+1} \ge \gamma_n$.

We first consider the case $n \le N$. Then the recursion formula
\eqref{eq:sigman} yields
\begin{equation} \label{eq:gamman}
\gamma_{n+1} \ge \gamma_n + \sum_{l=1}^{n-1} (b_l - \curv_{n-l})
\gamma_{n-l},
\end{equation}
and $\gamma_{n+1} \ge \gamma_n$ follows from Lemma \ref{lem:4facts}(d).

Finally, we consider the case $n > N$. If $N \ge 2$, the recursion
formula \eqref{eq:sigman}, Lemma \ref{lem:4facts}(a,b) and the
monotonicity of $\gamma_k$ yields
\begin{eqnarray*}
\gamma_{n+1} &\ge& \gamma_n + \left( \sum_{l=1}^{N-1} (b_l - \curv_{n-l})
\gamma_{n-l} \right) - \gamma_{n-N} \\
&\ge& \gamma_n + \sum_{l=1}^{N-2} (b_l - \curv_{n-l}) \gamma_{n-l}.
\end{eqnarray*}
Again, $\gamma_{n+1} \ge \gamma_n$ follows now from Lemma \ref{lem:4facts}(d).
If $N=1$ (i.e., $q =3$ or $q=4$), the recursion formula
\eqref{eq:sigman} simplifies considerably and, using Lemma \ref{lem:4facts}(c),
we conclude that
$$ \gamma_{n+1} \ge 2 \gamma_n - \gamma_{n-1} \ge \gamma_n, $$
finishing the proof of Theorem \ref{thm:growest}.
\end{proof}

\section{Proof of Theorem \ref{thm:treecomp}}
\label{sec:proofthm4}

\begin{proof}[Proof of Theorem \ref{thm:treecomp}]
We choose a vertex $v_0 \in \V$ and introduce the following functions
$m,M: \F \to \{ 0,1,2,\dots, \infty \}$:
\begin{eqnarray*}
m(f) &=& \min \{ d(w,v_0) \mid w \in \partial f \},\\
M(f) &=& \max \{ d(w,v_0) \mid w \in \partial f \}.
\end{eqnarray*}
Note that the face $f$ ``opens up'' at distance $m(f)$ and ``closes
up'' at distance $M(f)$ from $v_0$.  We call a face $f$ {\em finite},
if $M(f) < \infty$.

The idea of the proof is to ``open up'' successively every finite face
$f \in \F$ into an infinigon without violating the vertex bound. In
this way, we will build up a comparison tree $\T$ with the same vertex
bound $p$ and satisfying $\mu(\G) \le \mu(\T)$. It turns out, however,
that finite faces $f$ with more than one vertex in the sphere
$S_{M(f)}(v_0)$ cause problems in this ``opening up''
procedure. Therefore, we first modify the tessellation $\G$ by
removing all edges connecting two vertices $v,w$ at the same distance
to $v_0$.  The modified planar graph is denoted by $\G_0 =
(\V_0,\E_0,\F_0)$. The modification $\G \to \G_0$
is illustrated in Figure \ref{remvert}. (For convenience, the vertices
belonging to distance spheres $S_n(v_0)$ are arranged to lie on
concentric Euclidean circles around $v_0$.)

\begin{figure}[h]
  \begin{center}
    \psfrag{G}{\Large $\G$} \psfrag{G'}{\Large $\G_0$}
    \psfrag{v0}{$v_0$}
    \includegraphics[width=12cm]{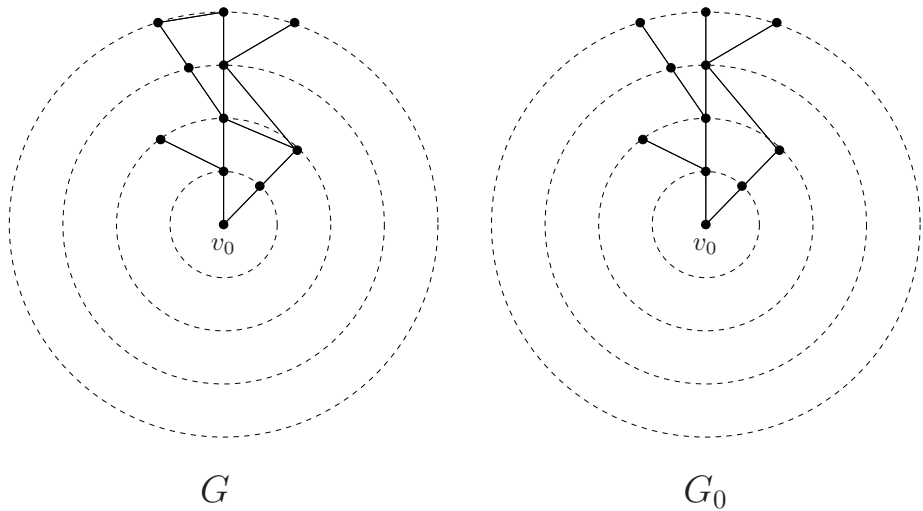}
  \end{center}
  \caption{Removing edges between vertices on the same spheres}
  \label{remvert}
\end{figure}

Note that none of the distance relations of the vertices in $\G_0$ to
the vertex $v_0$ are changed and that we still have $\Cut(v_0) =
\emptyset$. Moreover, the modified graph $G_0$ has a new set of faces
$\F_0$. Every finite face $f$ of $\G_0$ has now even degree, since $f$
opens up at a single vertex in the sphere $S_{m(f)}(v_0)$ and $f$
closes up at a single vertex in the sphere $S_{M(f)}(v_0)$.

We order all finite faces $f_0, f_1, f_2, \dots$ of $\G_0$ such that
we have
$$ M(f_0) \le M(f_1) \le M(f_2) \le ... $$
Next we explain the first step of our procedure, namely, how to open
up $f_0$ into an infinigon $\widetilde f_0$. Let $n = M(f_0) \ge 1$
and $w \in \partial f_0$ such that $d(w,v_0) = n$. Since $C(v_0) =
\emptyset$, we can find an infinite geodesic ray $w_0= w, w_1, w_2,
\dots \in \V$ such that $d(w_i,v_0) = n+i$. We may think of $v_0$ as
being the origin of the plane and of $w_0, w_1, \dots$ as being
arranged to lie on the positive vertical coordinate axis at heights
$n,n+1, \dots$ with straight edges between them. Now we cut our plane
along this geodesic ray, i.e., replace the ray by two parallel copies
of the ray and thus preventing the face $f_0$ from closing up at
distance $n$. In this way, $f_0$ becomes an infinigon, which we denote
by $\widetilde f_0$. The procedure
is illustrated in Figure \ref{openup}. Note that the vertices $w_i$
are replaced by two copies $w_i^{(1)}, w_i^{(2)}$, such that
$w_i^{(j)}$ is connected to $w_{i+1}^{(j)}$ for $j=1,2$ and
$w_i^{(1)}$ inherits all previous neighbors of $w_i$ at one side of
the ray and $w_i^{(2)}$ inherits all previous neighbors of $w_i$ at
the other side of the ray. In this way we obtain a new planar graph
$\G_1 = (\V_1,\E_1,\F_1)$.

\begin{figure}[h]
  \begin{center}
    \psfrag{G0}{\Large $\G_0$} \psfrag{G1}{\Large $\G_1$}
    \psfrag{v0}{$v_0$} \psfrag{f0}{$f_0$} \psfrag{F0}{$\widetilde
      f_0$}
    \includegraphics[width=12cm]{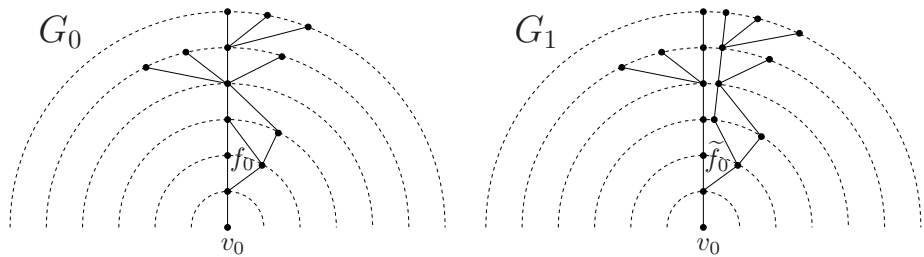}
  \end{center}
  \caption{Changing the finite face $f_0$ into an infinigon
    $\widetilde f_0$}
  \label{openup}
\end{figure}

The graph $\G_1$ is still connected and we obviously have
$\mu(\G_1,v_0) \ge \mu(\G_0,v_0) = \mu(\G,v_0)$. We then carry out the
same procedure with the face $f_1$ and the face $f_2$ and so on. The
limit is a connected tree $\T$ (since all faces of $\T$ are
infinigons) satisfying $\mu(\T,v_0) \ge \mu(\G,v_0)$ with vertex degrees
bounded by $p$. Of course, adding branches to make it a $p$-regular tree
only further increases the exponential growth. This finishes the proof
of Theorem \ref{thm:treecomp}. \end{proof}


\end{document}